\documentstyle[12pt,epsf]{article}

\def\reel{\hbox{{\rm R}\kern-1em\hbox{{\rm I} }}}
\def\relatif{\ \hbox{{\rm Z}\kern-.4em\hbox{\rm Z}}}
\def\nat{\hbox{{\rm N}\kern-1em\hbox{{\rm I} } }}
\def\comp{\hbox{{\rm C}\kern-.55em\hbox{{\rm I} } }}
\def\smallcomp{\hbox{\fiverm C}\kern-.35em{\hbox{\fiverm I}}}

\def\fudge{\mathchoice{}{}{\mkern.5mu}{\mkern.8mu}}
\def\bbc#1#2{{\rm \mkern#2mu\vbar\mkern-#2mu#1}}
\def\bbb#1{{\rm I\mkern-3.5mu #1}} \def\bba#1#2{{\rm #1\mkern-#2mu\fudge
#1}}
\def\bb#1{{\count4=`#1 \advance\count4by-64 \ifcase\count4\or\bba
A{11.5}\or \bbb B\or\bbc C{5}\or\bbb D\or\bbb E\or\bbb F \or\bbc
G{5}\or\bbb H\or \bbb I\or\bbc J{3}\or\bbb K\or\bbb L \or\bbb
M\or\bbb N\or\bbc O{5} \or \bbb P\or\bbc Q{5}\orrrr\b bb R\or\bbc
S{4.2}\or\bba T{10.5}\or\bbc U{5}\or    \bba V{12}\or\bba
W{16.5}\or\bba X{11}\or\bba Y{11.7}\or\bba Z{7.5}\fi}}
\def\rat{\hbox{{\rm Q}\kern-.70em\hbox{{\rm I} } }}

\newcommand{\be}{\begin{equation}}
\newcommand{\ee}{\end{equation}}

\newcommand{\ber}{\begin{eqnarray}}
\newcommand{\eer}{\end{eqnarray}}

\newcommand{\nin}{\noindent}

\def\bbb#1{{\rm I\mkern-3.5mu #1}} \def\bba#1#2{{\rm #1\mkern-#2mu\fudge
#1}}

\newcommand{\ein}{\eta\in{\cal{X}}_N}

\newcommand{\la}{\label}

\title{The simplest nearest-neighbor spin systems on regular graphs: Time
dynamics of the mean coverage function }

\author{{\bf Boris L. Granovsky}
\thanks{E-mail: mar18aa@techunix.technion.ac.il} \\
Department of Mathematics, Technion-Israel Institute of
Technology,\\ Haifa, 32000,Israel.}
\begin{document}
\maketitle \vskip 5cm

\nin American Mathematical Society 1991 subject classifications.

\nin Primary-60J27;secondary-60K35,82C22,82C26.

\nin Keywords and phrases: Interacting spin systems, Mean coverage
function, Voter models, Ergodicity, Spectral gap.

\newpage
\setcounter{equation}{0}
\begin{abstract}
\nin We establish  a characterization of the class of the simplest
nearest neighbor spin systems possessing the mean coverage
function (mcf) that obeys a  second order differential equation,
and derive explicit expressions for the mcf's of  the above
models. Based on these  expressions, the problem of ergodicity of
the  models is studied and bounds for their spectral gaps are
obtained.
\end{abstract}

\section{Introduction and Summary}

\nin It is commonly acknowledged that  even in the  case of a
simple infinitesimal interaction mechanism, a  description of the
transient behavior  of an interacting particle system (IPS)  is an
intractable mathematical problem in the theory of Markov
processes. In view of this, there is a continuing interest, both
in theory and applications, in seeking  solvable (in some sense)
models of IPS. One of the important  functionals of IPS is
undoubtedly, their  mean coverage function. In the present paper
we continue to study the behavior in time of the mean coverage
function  of a class of IPS called the simplest nearest neighbor
spin systems (SNNSS) on $s$-regular graphs. Namely, developing the
approach of Granovsky and Rozov \cite{GR}, we establish a
characterization of the class of SNNSS that posses a mean coverage
function satisfying a second order differential equation.
 \nin This is the main result of the present paper, stated in the
Theorem, Section 3. The theorem asserts that the above class
consists of the following four different modifications of the
basic voter model:
 noisy voter model,  noisy voter model with threshold $=2$ ( or
 $3$)
  on $2$- (resp.$3$) regular graphs,  a special case of a
general   threshold $=2$ model in one dimension and a degenerate
model with threshold=$s$ on $s$-regular graphs.

\nin It should be noted that  the first of these models is the
unique  SNNSS that has a mean coverage function satisfying a first
order differential equation. This was proven in \cite{GR}. In
Section 3 we  derive explicit expressions for the mean coverage
functions of the above  four models, by solving the corresponding
second order differential equations. The formulae obtained show
that adding a constant noise to flip rates results in considerable
change in transient behavior of the process. This matter is
discussed in Section 5.

\nin The next two sections are based on the aforementioned
formulae for the mean coverage function.
 Section 4 is devoted to the mean density function. We prove here
 that, when started from the product Bernoulli measure, the mean
density functions of the above processes do not depend on the size
of the graph. This remarkable property is used for the study of
ergodicity in the next section. Section 4 contains also a
historical sketch of research related to the subject.

\nin It is clear that transient behavior of the mean coverage
function, which is of interest in itself,  also provides
information on the long-time properties of the process considered.
In view of this, the last Section 5 is devoted to ergodicity and
bounding the spectral gap for the class of models defined in
Theorem. We give here a positive answer on the  open problem about
ergodicity of threshold $=3$ noisy voter model on $3$-regular
graphs, for some values of parameters.  Based on the expressions
derived in Section 3, we obtain the upper bounds for the spectral
gap of the four SNNSS. These  upper bounds are compared with the
lower bounds given by the $\epsilon- M>0$ condition (for
references see \cite{Lig1}, p.31).

\nin Finally, note that in the  course of the proof of the Theorem
we derived identities that hold for a  coverage  of sites of a
regular graph by $0'$s and $1'$s. These identities might be
helpful in the study of other problems related to time dynamics of
SNNSS.

 \nin Most of the  notation and language of our paper
have been  adopted from the seminal monograph on IPS \cite{Lig1},
by Liggett.

\section{Background}

\nin   We consider throughout the paper a SNNSS on a $s$-regular
graph $G$ of finite size $N,$ with the set of vertices(sites)
$V=\{x\}.$ Recall that a graph is called $s$-regular if each of
its vertices has $s$ neighbors. By SNNSS
 we mean a time homogeneous Markov process $\varphi_t, \quad t\ge 0$
 with state
space ${\cal{X}}_N$ = $\{ 0,1\}^V=\{\eta\}$ and the  infinitesimal
time dynamics given by (\ref{00}) below. The elements
$\eta=\{\eta(x), \quad x\in V\}$ of  ${\cal{X}}_N$ are called
configurations. We will say that a site $x\in V$ is
occupied(resp., empty) in the configuration $\ein, $ if $\eta(x)$
is $1$ (resp., 0). The SNNSS are featured by the property that the
flip rate $c(x,\eta)$ of a spin at a site $x\in V$ in a
configuration $\ein $ depends only on the number $k(x,\eta)$ of
occupied neighbors of $x$ in the configuration $\eta.$ Formally,
\be
c(x,\eta)= \lambda_k(1-\eta(x))+\mu_k\eta(x), \quad
k=k(x,\eta),\quad x\in V, \quad \ein, \la{01} \ee

\nin where  $\lambda_k, \quad k=0,1,\dots, s$ (resp., $\mu_k,
\quad k=0,1,\dots, s$ ) are the rates of the infinitesimal
transitions $ 0\rightarrow 1$ (resp., $1\rightarrow 0$) at a given
site in a given configuration. Finally, denoting by $\eta_x$ the
configuration obtained from $\eta$ by flipping the spin at the
site $x,$ the above assumptions conform to the following
infinitesimal time dynamics of $\varphi_t, \quad t\ge 0:$
\be
Pr(\varphi_{t+\Delta t} =\eta_x \mid \varphi_{t}=\eta)=
c(x,\eta)\Delta t + o(\Delta t), \ \ \quad \Delta t\ge 0, \quad
t\ge 0, \quad x\in V,\quad \ein, \la{00}
 \ee
\nin where $\frac{o(\Delta t)}{\Delta t}\to 0$, as $\Delta t\to
0.$

\nin So, to compare with a variety of the so called biased models,
( see e.g. Madras, Schinazi and Schonmann \cite{MS}) SNNSS is a
spatially homogeneous process.

\nin It is known that in the above setting the process $\varphi_t,
\quad t\ge 0$ is fully defined by the $2s+2$ parameters
$\lambda_k\ge 0 , \mu_k\ge 0, \quad k=0,\ldots, s.$ Namely, the
generator $\Omega$ of the process is given by

\be
\Omega f(\eta)=\sum_{x\in V} c(x,\eta)(f(\eta_x)-f(\eta)),\quad
f\in C(\cal{X}_N), \quad \ein, \la{02} \ee

\nin where $C({\cal{X}}_N)$ is the class of bounded functions
$f:{\cal{X}}_N\rightarrow R.$

\nin Denote $\varphi^{(\eta)}_t, \quad t\ge 0,\  \ein $ the SNNSS
starting from a configuration $\eta$ and

\nin $M_f^{(\eta)}(t)=Ef(\varphi^{(\eta)}_t), \ t\ge 0, \ f\in
C({\cal{X}}_N),\ \ein.$ As in Granovsky and Rozov \cite{GR}, our
starting point will be the following assertion that is a
straightforward consequence of the Hille- Yosida theorem.

\nin {\bf Proposition 1.} The function $M_f^{(\eta)}(t),\  t\ge 0$
satisfies, for all $\ein,$ a
 linear
differential equation of order $l,\ ( l\ge 1)$
\be
\frac{d^lM_f^{(\eta)}(t)}{dt^l}= \sum_{i=0}^{l-1}
A_i\frac{d^iM_f^{(\eta)}(t)}{dt^i}+ B,\quad t\ge 0,\quad \ein
\la{03} \ee

\nin with coefficients $A_i,\quad i=0,\ldots, l-1$ and $ B$ that
do not depend on $\ein$ and $t\ge 0,$  iff the generator $\Omega$
of the Markov process considered obeys the condition

\be
\Omega^l f=\sum_{i=0}^{l-1} A_i\Omega^i f+B, \quad \ein, \la{04}
\ee

\nin where $\Omega^{i+1}:=\Omega(\Omega^i), \ i=0,1,\ldots.$ \nin

\nin Our subsequent study of the characterization problem
described in the previous section is  based on the fact that
(\ref{03}) is equivalent to (\ref{04}).

\section{Main result}

\nin  The coverage of the graph $G$ by a configuration $\ein$ is
the function $\vert\eta\vert:{\cal{X}}_N\rightarrow
R^+,\vert\eta\vert=\sum_{x\in V}\eta(x)$ and
$M^{(\eta)}(t):=E\vert \varphi^{(\eta)}_t\vert, \ t\ge 0, \ein $
is called  the mean coverage function of the process
$\varphi_t^{(\eta)}, \quad t\ge 0.$  The function $M^{(\eta)}(t),
t\ge 0$ is one of the most important functionals in applications.
 In Granovsky, Rolski, Woyczinski and Mann \cite{GW} and Belitsky,
Granovsky \cite{BG} the function was studied in the context of
adsorption - desorption process given by

\nin $\lambda_k=\lambda>0, \mu_k\ge 0, \ k=0,\ldots,s.$  It was
observed there that the function $M^{(\eta)}(t), t\ge 0,$ has a
saddle point, under certain conditions on parameters of the
process.

\nin Our main objective will be to describe the class of SNNSS
satisfying (\ref{03}) with $f=\vert\eta\vert$ and $l=2.$ For $l=1$
the problem was posed and solved in \cite{GR}.

\nin We introduce some more notation. Denote
\be
g_i(\eta):=\Omega^i(\vert\eta\vert), \ i=0,1,\ldots, \ \ein
\la{ab}\ee

\nin  to obtain from (\ref{02})

\be
g_1(\eta)= \sum_{x\in V} c(x,\eta)(1-2\eta(x)),\quad \ein. \la{05}
\ee

\nin In view of our objective, we will need to unlock the
structure  of $g_2$ .

\nin Let $ D$ be a nonempty subset of $V.$ We will say that $y\in
V$ is a neighbor of $D:$ $y\sim D$, if $y\notin D$ and  $y$ is a
neighbor of at least one site in $D,$ and we denote $\delta_1( D)$
the set of all neighbors of the subset $D.$ In particular, by
$\delta_1(x)$ we denote the neighborhood of $x\in V.$ We  also
define $\delta_i(x)=\delta_1(\delta_{i-1}(x)), \ \ i=1,2, \ldots,
\ \delta_0(x)=\{x\},  \  x\in V.$

 \nin Next, for any $x\in V$ define the difference operator $\Delta_x
:C({\cal{X}}_N)\to C({\cal{X}}_N)$

\be
\Delta_x f(\eta)=f(\eta_x)-f(\eta), \quad f\in C({\cal{X}}_N),
\quad x\in V \ee

\nin and write  $\Delta^{(2)}_{x,y} f=\Delta_x\Delta_y f, \quad
f\in C({\cal{X}}_N),\quad x,y\in V.$ Then, by our  definition
(\ref{ab}) and (\ref{02}) we have

\be g_2(\eta)=\sum_{x\in V}c(x,\eta)\Delta_xg_1(\eta) ,\quad \ein.
\la{e} \ee

\nin Further, it follows from (\ref{05}) and (\ref{01}) that

\begin{eqnarray}
\Delta_y g_1(\eta)=\sum_{x\in \delta_1(y)} (1-2\eta(x))\Delta_y
c(x,\eta)-\nonumber\\
\Big(c(y,\eta_y)+c(y,\eta)\Big)(1-2\eta(y)\big),\quad \ein,
\label{1}
 \end{eqnarray}

\nin  for any $y\in V.$ Since $\Delta_{x,y}=\Delta_{y,x}, \ x,y\in
V,$  (\ref{1}) implies the important fact that

\be
\Delta^{(2)}_{x,y} g_1(\eta)=0,\quad \ein , \label{3} \ee

\nin whenever $ x\not\in \delta_1(y)\bigcup \delta_2(y)$ and
$x\neq y.$ This and (\ref{e}) give

\begin{eqnarray}\
 \Delta_y g_2(\eta)=
 \Delta_y\Big[c(y,\eta)\Delta_yg_1(\eta)\Big]+\nonumber\\
  \sum_{x\in \delta_1(y)\bigcup
\delta_2(y)}\Delta_y \Big[c(x,\eta)\Delta_x g_1(\eta)\Big], \quad
y\in V, \quad \ein. \la{2} \end{eqnarray}

\nin For the proof of our main result, stated in the Theorem in
the sequel, we need to impose the following two conditions  on
$s$-regular graphs $G$ considered.

\nin (i.) First, we assume that $G$ is triangular free graph,
which means that if $x,y,z\in V: y,z\sim x,$ then $ y,z$ are not
neighbors. The second condition is a technical one.

\nin (ii.) We assume the existence of a pair of vertices $y,z\in
V$ s.t. $ z\in \delta_3( y)$ and the two sets of vertices
$E_{1,2}:= \delta_1(y)\cap \delta_2(z)$ and
$E_{2,1}:=\delta_2(y)\cap\delta_1(z)$ are singletons.

\nin Observe that the conditions (i) and (ii) are satisfied e.g.,
when $G$ is an $s$-regular tree or $G=Z^d, \quad d\ge 1.$

\nin For the purpose of establishing our characterization result
we  employ a  technique that is presented below. We start with the
notations adopted from \cite{BG} and \cite{GR} . Denote
$n_k^{(i)}= n_k^{(i)}(\eta), \ k=0,\ldots,s, \ i=0,1$ the number
of occupied (i=1) (resp., empty (i=0)) sites having $k$ occupied
neighbors in a configuration $\ein,$ and let
${\bf{n}}^{(0)}_k={\bf{n}}^{(0)}_k(\eta), {\bf
n}^{(1)}_k={\bf{n}}^{(1)}_k(\eta), \quad k=0,\ldots, s$ denote the
corresponding sets of vertices $x\in V$ in a configuration $\ein.$
Finally, we denote $V^{(i)}=V^{(i)}(\eta), i=0,1$ the set of all
empty (resp. occupied) sites in $\ein.$

\nin Then $g_1$ defined  by (\ref{05}) can be expressed as

\be
g_1(\eta)= \sum_{k=0}^s (\lambda_kn_k^{(0)}- \mu_kn_k^{(1)}),
\quad \ein. \la{97} \ee

\nin The following  identities that are valid for any $s$-regular
graph  will be crucial for our subsequent study:
\be
P=P(\eta):=\sum_{x\in\bf{n_s^{(0)}}} \Delta_x
n_0^{(1)}=\sum_{x\in\bf{n_0^{(1)}}} \Delta_x n_s^{(0)}, \quad \ein
\label{38} \ee

\be
Q_0=Q_0(\eta):=\sum_{x\in V^{(0)}}\Delta_x n_s^{(0)}= -n_s^{(0)}+
n_{s-1}^{(0)},\quad \ein \la{39}
 \ee

\be
Q_1=Q_1(\eta):=\sum_{x\in V^{(1)}}\Delta_x
n_s^{(0)}=-sn_s^{(0)}+n_s^{(1)},\quad \ein \la{40}\ee

\be
\sum_{x\in\bf{n_s^{(0)}}} \Delta_x n_s^{(0)}=-n_s^{(0)},\quad
\sum_{x\in\bf{n_0^{(1)}}} \Delta_x n_0^{(1)}=-n_0^{(1)}, \quad
\ein. \label{100} \ee

\nin The proof of the identities can be obtained after some
thought  from the preceding definitions. We also write

\be
R_i=R_i(\eta):=\sum_{x\in V^{(i)}}\Delta_x n_0^{(1)},\quad i=0,1,
\quad\ein. \la{R} \ee

\nin Let  $\bar{\eta}$ be the configuration obtained by flipping
the spins at all sites $x\in V$ in a configuration $\ein.$ Then we
have $n_k^{(i)}(\bar{\eta})=n_{s-k}^{(1-i)}(\eta),\quad
k=0,1\ldots,s,\quad i=0,1, \quad \ein$ and, consequently,

\be
R_i(\eta)=Q_{1-i}(\bar{\eta}),\quad \ein, \quad i=0,1. \la{41} \ee

\nin Now we are in a position to state the following

 \nin {\bf
Lemma.} The identity (in $\eta$)

\be
n_s^{(1)}+n_{s-1}^{(0)} -n_1^{(1)}-n_0^{(0)}=
(n_s^{(0)}-n_0^{(1)})F_1+ \vert\eta\vert F_2+ F_3,\quad \ein,
\la{80} \ee

\nin where $F_i, \ i=1,2,3$ are coefficients that do not depend on
$\eta,$ holds iff $s=2,3.$ In  both cases of $s$,\
$F_2=2,\ F_3=-N,$ while 
\be
F_1=\left \{
\begin{array}{ll}
-3, & {\rm if~} s=2 \cr -2, & {\rm if~} s=3.
\end{array}
\right. \la{89}\ee

\nin {\bf Proof.} We put in (\ref{80}) first $\eta=\emptyset$ and
then $\eta=\bar{\emptyset}$ to  find $F_2$ and $ F_3.$ Now
consider the case $s\ge 3.$  Due to the fact that the graph
considered is triangular free, we have for any $x\sim y, \ x,y\in
V,$

\be
n_s^{(1)}(\emptyset_{x,y})= n_{s-1}^{(0)}(\emptyset_{x,y})=
n_0^{(1)}(\emptyset_{x,y})=n_{s}^{(0)}(\emptyset_{x,y})=0, \quad
n_0^{(0)}(\emptyset_{x,y})=N-2s,\quad
n_1^{(1)}(\emptyset_{x,y})=2. \la{811}\ee

\nin Substituting this in (\ref{80}), gives $-2-(N-2s)= 4-N,$
which says that (\ref{80}) does not hold for $s>3.$

\nin Further, if $s=3 $ and $\eta=\emptyset_{x}, \ x\in V,$ then

\be
n_3^{(1)}(\emptyset_{x})= n_{2}^{(0)}(\emptyset_x)=
n_1^{(1)}(\emptyset_x)=n_{3}^{(0)}(\emptyset_x)=0, \
n_0^{(0)}(\emptyset_x)=N-4,\ n_0^{(1)}(\emptyset_x)=1. \la{81}\ee

\nin  The latter  implies $F_1=-2.$

\nin In the case $s=2$ we have $n_{s-1}^{(0)}(\emptyset_{x,y})=2,$
and  the same argument as before gives $F_1=-3.$ Finally, it is
left to show that the identity (\ref{80}) indeed holds for
$s=2,3.$ We use the relationship $$ 2\vert\eta\vert -N=
\vert\eta\vert -(N-\vert\eta\vert)= \sum_{k=0}^s (
n_k^{(1)}-n_k^{(0)})$$ to obtain for $s=3$
\be
n_3^{(1)}+n_{2}^{(0)} -n_1^{(1)}-n_0^{(0)}+
2(n_3^{(0)}-n_0^{(1)})-2\vert\eta\vert+N=\sum_{k=0}^3(kn_k^{(1)}-
(s-k)n_k^{(1)})=0,\quad \ein \la{82}, \ee

\nin where the last equation follows from the identity
$\sum_{k=1}^sk(n_k^{(0)}+n_k^{(1)})=s\vert\eta\vert, \ \ein$ that
is valid for all $s$-regular graphs. The same argument proves the
assertion  for $s=2.$ $\clubsuit$ \vskip .5cm

\nin Finally, we will distinguish the following modifications of
the Basic Voter model :

\nin {\bf Noisy Voter Model.} $\lambda_k-\lambda_{k-1}=
\mu_{k-1}-\mu_k=d, \ k=1,\ldots, s.$ The model was introduced in
\cite{GR} and intensively studied in \cite{GM}. Here the noise is
given by the two parameters $h_1=\lambda_0,\ h_2=\mu_0-sd$ added
to the basic voter model (see \cite{Lig1},\cite{Lig2}):
$\lambda_k=kd, \ \mu_k=(s-k)d, \ k=0,1, \ldots, s$

\nin  Note that in \cite{Lig1} Ex.2.5, p.136, it is considered a
general (i.e. not necessarily the nearest neighbor) version of
voter model with noise.

\nin {\bf Noisy  Voter Model with Threshold} $=q\ (1\le q\le s).$

\be
\lambda_k=\mu_{k+s-q+1}=h\ge 0,\ k=0,\ldots, q-1, \quad
\lambda_k=\mu_{k-q}= h+a\ge 0, \ k=q,\ldots, s.  \ee

\nin This is  the simplest case of a nonlinear voter model. In the
case $h=0$ (the absence of noise),  the model was suggested by Cox
and Durrett in \cite{CD}. ( For updated references see
\cite{Lig2}). In \cite{CD} it was also considered the threshold
voter model with noise added to the death rates only. If $q=s,$
then by scaling all the rates by the factor $(2h+a)^{-1}$  the
model becomes
 the nearest neighbor Majority Vote Process (\cite{ Lig1},
Ex.4.3(e),p.33 and Ex. 2.12,p. 140).

 \nin
{\bf Generalized Threshold Model with threshold $=q \ (1\le q\le
s)$.}  The model is obtained from the previous one by adding a
constant either  to  $s-q+1$ birth rates $\lambda_k, \quad
k=q,\ldots,s,$ or to $s-q+1$ death rates $\mu_k, \quad k=0,\ldots,
q-s.$ Explicitly,

$$ \lambda_k=\mu_{k+s-q+1}=h\ge 0,\ k=0,\ldots, q-1, \quad
\lambda_k= h+a\ge 0, \ k=q,\ldots, s,$$ \be \mu_{k}= h+b\ge 0, \
k=0,\ldots, s-q. \ee
\newpage
\nin {\bf Theorem}

\nin  The mean coverage function $M^{(\eta)}(t),\  t\ge 0$ of a
SNNSS $\varphi_t,\ t\ge 0$ satisfies, for all $\ein,$ a second
order linear differential equation

\be
\frac{d^2M^{(\eta)}(t)}{dt^2}= A_1\frac{dM^{(\eta)}(t)}{dt}+A_0
M^{(\eta)}(t)+B,\quad t\ge 0,\quad \ein \la{cd} \ee

\nin with coefficients $A_0,A_1, B$ that do not depend on $\ein$
and $t\ge 0,$ iff $\varphi_t,\ t\ge 0$ is  one  of the following
four models $(\mathbf{C_1}) -(\mathbf{C_4})$ :

\nin $(\mathbf{C_1})$ A noisy voter model.

\nin $(\mathbf{C_2})$  A generalized threshold model with
threshold $=s$ and $h= ab=0$ or $h=0,\ a=b.$

\nin $(\mathbf{C_3})$ A threshold noisy voter model with threshold
$=s,$ when $s=2,3$

\nin $(\mathbf{C_4})$  A generalized threshold model with
threshold $s=2$ and $h, a,b: h(a+b)=ab,
 \ h\ge 0,\ h+a\ge 0, \ h+b\ge 0.$

\nin {\bf Proof.} By virtue of Proposition 1, (\ref{cd}) is
equivalent to
\be
g_2(\eta)= A_1 g_1(\eta) +  A_0\vert\eta\vert + B, \quad \ein,
\la{90}\ee

\nin The main difficulty is to prove  that (\ref{90}) implies one
of the four conditions $(\mathbf{C_1}) -(\mathbf{C_4})$ on the
rates of $\varphi_t,\ t\ge 0$

 \nin If (\ref{90}) holds, then, by
(\ref{3})

\be
\Delta^{(2)}_{y,z} g_2(\eta)= A_1\Delta_{y,z}^{(2)} g_1(\eta)=0,
\quad y\in V,\quad z\in \delta_3(y), \quad \ein. \la{91}
\ee

\nin From the other hand, we get from (\ref{2}) and (\ref{3})
 \be
\Delta^{(2)}_{y,z} g_2(\eta)= \sum_{x\in \delta_1(y)\bigcup
\delta_2(y)}\Delta_{y,z}^{(2)} \Big[c(x,\eta)\Delta_x
g_1(\eta)\Big], \quad y\in V, \quad z\in \delta_3(y), \quad \ein.
\ee

\nin In view of  (\ref{3}) this  gives

\be \Delta^{(2)}_{y,z} g_2(\eta)= \sum_{x\in E_{1,2}}
\Delta_y\Big[c(x,\eta)
  \Delta^{(2)}_{z,x} g_1(\eta)\Big]
+ \sum_{x\in E_{2,1}}\Delta_z\Big[c(x,
 \eta)\Delta^{(2)}_{x,y} g_1(\eta)\Big], \\
 \quad z\in \delta_3(y),\quad
 \ein,
\label{4} \ee

\nin  where we denoted $E_{1,2}:= \delta_1(y)\cap \delta_2(z)$ and
$E_{2,1}:=\delta_2(y)\cap\delta_1(z).$ We also derive from
(\ref{1})

\be
\Delta^{(2)}_{x,z}
g_1(\eta)=\sum_{u\in\delta_1(x)\cap\delta_1(z)}(1-2\eta(u))
\Delta^{(2)}_{x,z} c(u,\eta), \quad x\in \delta_2(z),\quad \ein
 \ee

\nin and

\be
\Delta^{(2)}_{x,y} g_1(\eta)=
\sum_{u\in\delta_1(x)\cap\delta_1(y)}(1-2\eta(u))\Delta_{
x,y}^{(2)}c(u,\eta), \quad x\in \delta_2(y), \quad \ein. \label{6}
\ee

\nin We substitute now these expressions in (\ref{4}) to obtain

\begin{eqnarray}
\Delta^{(2)}_{y,z} g_2(\eta)= \sum_{x\in E_{1,2}}
 \Delta_y\Big[ c(x,\eta)\sum_{u\in\delta_1(x)\cap\delta_1(z)}
 (1-2\eta(u))
 \Delta_{x,z}^{(2)}c(u,\eta)\Big]+\nonumber\\ \sum_{x\in E_{2,1}}
 \Delta_z\Big[ c(x,\eta)
 \sum_{u\in\delta_1(x)\cap\delta_1(y)}
 (1-2\eta(u))\Delta_{
x,y}^{(2)}c(u,\eta)\Big],\quad y\in V, \quad z\in \delta_3(y),
\quad \ein. \label{7}
\end{eqnarray}

\nin Our immediate aim is to find conditions on the parameters of
a SNNSS,  imposed by the requirement

\be
\Delta^{(2)}_{y,z} g_2(\eta)=0, \quad y\in V, \quad z\in
\delta_3(y), \quad \ein.\la{92}  \ee

\nin Let, in accordance with the  assumption (ii)  on $G,$ the
vertices $y,z$ in (\ref{7}) be such that

\be E_{1,2}=\{u_1\}, \quad E_{2,1}=\{u_2\}, \la{101}\ee

\nin where $u_1, u_2\in V.$ Then, in view of the above definition
of the vertices $u_1, u_2,$ (\ref{7}) becomes

\begin{eqnarray}
\Delta^{(2)}_{y,z} g_2(\eta)&=&
 \Big(\Delta_y c(u_1,\eta)\Big)
 (1-2\eta(u_2))
 \Delta^{(2)}_{u_1,z}c(u_2,\eta) \nonumber\\
 &+&\Big(\Delta_z c(u_2,\eta)\Big)
 (1-2\eta(u_1))\Delta_{
u_2,y}^{(2)}c(u_1,\eta), \quad \ein, \quad z\in
\delta_3(y).\label{8}
\end{eqnarray}

\nin The last expression will be our main tool in  the subsequent
study.

\nin   We see from (\ref{8}) that $$\Delta^{(2)}_{y,z}
g_2(\eta_y)=-\Delta^{(2)}_{y,z} g_2(\eta), \quad
\Delta^{(2)}_{y,z} g_2(\eta_z)=-\Delta^{(2)}_{y,z} g_2(\eta),
\quad z\in \delta_3(y),\quad \ein.$$ In view of this we set in
(\ref{8}), $\eta(y)=\eta(z)=0$ . We also agree to write $\Delta
(\bullet)_k=(\bullet)_{k+1}-(\bullet)_k$ and $\Delta^{(2)}
(\bullet)_k=(\bullet)_{k+2}-2(\bullet)_{k+1} + (\bullet)_{k},$
where $(\bullet)$ is either $\lambda$ or $\mu.$

\nin Now we will be attempting to find the explicit form of
 (\ref{8}) in the following three cases of $\ein $ that exhaust all the
possibilities. For brevity, we denote $k_i=k(u_i, \eta), \ \ein. $
It is important to note that $\delta_1(u_1)\bigcap \delta_1(u_2)$
is the empty set, since $G$ is triangular free.

\nin {\bf Case 1.} $\ein:\eta(u_1)=\eta(u_2)=0.$

\be
 \Delta^{(2)}_{y,z}g_2(\eta)=
\Big(\Delta\lambda_{k_1}\Big)\Big(\Delta^{(2)}\lambda_{k_2}\Big)+
 \Big(\Delta\lambda_{k_2}\Big)\Big(\Delta^{(2)}\lambda_{k_1}\Big),
\quad z\in \delta_3(y),\quad \ein, \label{9} \ee

\nin where  $0\le k_1,k_2\le s-2$.

\nin {\bf Case  2.} $\ein:\eta(u_1)=\eta(u_2)=1.$

\be
\Delta^{(2)}_{y,z}g_2(\eta)=
\Big(\Delta\mu_{k_1}\Big)\Big(\Delta^{(2)}\mu_{k_2-1}\Big)+
 \Big(\Delta\mu_{k_2}\Big)\Big(\Delta^{(2)}\mu_{k_1-1}\Big),
\quad z\in \delta_3(y),\quad \ein, \label{10} \ee

\nin where $1\le k_1,k_2\le s-1.$

\nin {\bf Case 3.} $\ein:\eta(u_1)=0, \ \ \eta(u_2)=1.$

\be
 \Delta^{(2)}_{y,z}g_2(\eta)=
\Big(\Delta\lambda_{k_1}\Big)\Big(\Delta^{(2)}\mu_{k_2}\Big)+
 \Big(\Delta\mu_{k_2}\Big)\Big(\Delta^{(2)}\lambda_{k_1-1}\Big),
\label{11} \ee

\nin where $1\le k_1\le s-1,\quad  0\le k_2\le s-2.$

\nin We now know from  the  three cases considered,  that the
condition (\ref{92}) implies


\ber
 \Big(\Delta\lambda_{k_1}\Big)\Big(\Delta^{(2)}\lambda_{k_2}\Big)+
 \Big(\Delta\lambda_{k_2}\Big)\Big(\Delta^{(2)}\lambda_{k_1}\Big)&=&
 0,\quad 0\le k_1,k_2\le s-2 \label{14a}\\
 \Big(\Delta\mu_{k_1}\Big)\Big(\Delta^{(2)}\mu_{k_2-1}\Big)+
 \Big(\Delta\mu_{k_2}\Big)\Big(\Delta^{(2)}\mu_{k_1-1}\Big)&=& 0,
\quad 1\le k_1,k_2\le s-1 \label{14b}\\
\Big(\Delta\lambda_{k_1}\Big)\Big(\Delta^{(2)}\mu_{k_2}\Big)+
 \Big(\Delta\mu_{k_2})\Big(\Delta^{(2)}\lambda_{k_1-1}\Big)&=&
 0,\quad 1\le k_1\le s-1,\nonumber\\
 & & \ \ \ \ \ \ 0\le k_2\le s-2.
 \label{14}
\eer

\nin Setting in (\ref{14a}), (\ref{14b}) $k_1=k_2=k$ gives \ber
\Big(\Delta\lambda_{k}\Big)\Big(\Delta^{(2)}\lambda_{k}\Big)&=&
0,\quad 0\le k\le s-2 \label{15a]}\\
\Big(\Delta\mu_{k}\Big)\Big(\Delta^{(2)}\mu_{k-1}\Big)&=& 0, \quad
1\le k\le s-1. \label{15} \eer

\nin Using this fact,  we multiply  the equations
(\ref{14a}),(\ref{14b})  by $\Delta\lambda_{k_1}$ and  by
$\Delta\mu_{k_1}$ correspondingly, to obtain

\ber
 \Big(\Delta\lambda_{k_1}\Big)\Big(\Delta^{(2)}\lambda_{k_2}\Big)&=&
 0, \quad 0\le k_1,k_2\le s-2  \nonumber\\
\Big(\Delta\mu_{k_1}\Big)\Big(\Delta^{(2)}\mu_{k_2-1}\Big)&=&
0,\quad 1\le k_1,k_2\le s-1.
 \label{16}
\eer

\nin Since, by our definition,
$\Delta^{(2)}(\bullet)_{k_2}=\Delta(\bullet)_{k_2+1}-\Delta(\bullet)_{k_2},$
  the equations (\ref{16}) imply for $s>2$

\ber \Delta^{(2)} \lambda_k&=&0, \quad 0\le k\le s-3. \nonumber\\
\Delta^{(2)} \mu_k&=&0, \quad 1\le k\le s-2. \label{18} \eer

\nin Thus, we have for $s>2$ \ber \Delta\lambda_k:&=
&d_\lambda,\quad k=0,\ldots,s-2\nonumber\\ \Delta\mu_k:&=& d_\mu,
\quad k=1,\ldots,s-1\label{19} \eer

\nin Finally, in view of  (\ref{18}) and (\ref{19}), we obtain
from (\ref{14a})-(\ref{14}) for all $s\ge2$

\ber d_\lambda\Delta^{(2)}\lambda_{s-2}&=& 0 \nonumber\\ d_\mu
\Delta^{(2)}\mu_{0}&=& 0\nonumber\\
\Big(\Delta\lambda_{s-1}\Big)\Big(\Delta^{(2)}\mu_{0}+
 \Big(\Delta\mu_{0}\Big)\Big(\Delta^{(2)}\lambda_{s-2}\Big)&=&
 0 .
\label{20} \eer

\nin Summarizing the preceding argument we conclude that
(\ref{19}) together with (\ref{20})  are necessary and sufficient
for (\ref{92}). Our next step will be devoted to show that the
conditions  (\ref{19}), (\ref{20}) on the parameters $\lambda_k,\
\mu_k, \ k=0,\ldots, s$ imply  one of the conditions
$(\mathbf{C_1}) -(\mathbf{C_4}).$

\nin Assume first that in (\ref{20}) $d_\lambda\neq 0.$ Then we
should have  $\Delta^{(2)}\lambda_{s-2}=\Delta
\lambda_{s-1}-\Delta\lambda_{s-2}=0,$ and, consequently, in view
of (\ref{19}), $\Delta \lambda_{s-1}=d_\lambda\neq 0$. Hence, in
view of the last equation in (\ref{20}), we obtain
\be
\Delta^{(2)}\mu_{0}=\Delta^{(2)}\lambda_{s-2}=0. \label{21}
 \ee

\nin By the same argument, (\ref{21}) should also hold under the
assumption $d_\mu\neq 0.$ (\ref{21}) together with (\ref{19}) is
equivalent to saying that the flip rates $\lambda_k, \quad
k=0,\ldots, s$ and  $\mu_k, \quad k=0,\ldots, s$ form arithmetical
progressions:

\ber \Delta\lambda_{k}&=& d_\lambda, \quad k=0,\ldots,
s-1\nonumber\\
 \Delta\mu_{k}&=&d_\mu, \quad k=0,\ldots, s-1.
 \label{22}
\eer

\nin So, assuming $d_\lambda\neq 0,$ one has
\be
 \Delta_y c(x,\eta)= (1-2\eta(y))\Big(d_\lambda(1-\eta(x))+
 d_\mu\eta(x)\Big), \quad x\in \delta_1(y), \quad x,y\in V,\quad
 \ein.
 \label{24}
 \ee
\nin Now (\ref{24}) and  (\ref{1}) imply

\be
 \Delta_y g_1(\eta)= (1-2\eta(y))\Big(sd_\lambda
 -2k(y,\eta)(d_\mu+
 d_\lambda)- \lambda_0-\mu _0\Big),\quad y\in V, \quad \ein,
\label{23} \ee

\be
\Delta^{(2)}_{y_2,y}c(x,\eta)=0,\quad y_2\in \delta_2(y),\quad
x\in V,\quad \ein,
 \ee

\nin and

 \be
 \Delta^{(2)}_{y_2,y} g_1(\eta)=0, \quad y_2\in \delta_2(y),\quad
\ein. \la{95}
 \ee

\nin Hence, by virtue of (\ref{90}) it follows from (\ref{95})
that \be  \Delta^{(2)}_{y_2,y} g_2(\eta)=0, \quad y_2\in
\delta_2(y),\quad \ein.
 \ee

\nin By (\ref{22}) the latter is equivalent to

$$ \Delta^{(2)}_{y_2,y}\Big(c(y_1,\eta)\Delta_{y_1}g_1(\eta)\Big)
=-4(d_\lambda +d_\mu)
(1-2\eta(y))(1-2\eta(y_2))\Big((1-\eta(y_1))d_\lambda-
\eta(y_1)d_\mu\Big)=0,$$
 \be y_2\in \delta_2(y),  \quad y_1\in
\delta_1(y), \quad \ein. \label{26}
 \ee

\nin This implies $d_\lambda+d_\mu=0,$ which by (\ref{22}),
 corresponds to the noisy voter model ($\mathbf{C_1}$).
 Since the same conclusion is valid under the assumption $d_\mu\neq 0,$
 it is left to assume that
$d_\lambda=d_\mu=0.$ In this case it follows from (\ref{19}) that
the parameters of  the  SNNSS are of the form

 \be
\lambda_k=\lambda, \quad k=0,\ldots, s-1, \quad \mu_k=\mu, \quad
k=1,\ldots,s, \quad \lambda_s =\lambda+a, \quad \mu_0=\mu+b,
\la{96}\ee

\nin where $a,b\in R$ are such that $\lambda+a\ge 0, \quad
\mu+b\ge 0.$

\nin Note that in the case considered all three  conditions
(\ref{20}) are satisfied, because (\ref{96}) implies
$\Delta^{(2)}\lambda_{s-2}=a, \ \ \Delta^{(2)}\mu_{0}=-b.$
  In view of the relationships $ \sum_{k=0}^s
n_k^{(0)}=N-\vert\eta\vert $  and

\nin $\sum_{k=0}^s n_k^{(1)}=\vert\eta\vert,$  (\ref{97}) yields
for the model (\ref{96})

\be
g_1(\eta)=\lambda N- (\lambda+\mu)\vert\eta\vert + an_s^{(0)}-
bn_0^{(1)}, \quad \ein. \label{36}
 \ee

\nin By (\ref{02}) we also have

 \ber \Omega(n_k^{(i)})=\lambda\sum_{x\in V^{(0)}\backslash
\bf{n_s^{(0)}}} \Delta_x n_k^{(i)} + \lambda_s
\sum_{x\in\bf{n_s^{(0)}}} \Delta_x n_k^{(i)}+\mu\sum_{x\in
V^{(1)}\backslash \bf{n_0^{(1)}}} \Delta_x n_k^{(i)}&+&
\nonumber\\
 \mu_0\sum_{x\in\bf{n_0^{(1)}}} \Delta_x n_k^{(i)}, \quad \ein,
\quad k=0,\ldots ,s,\quad i=0,1. \label{37} \eer

\nin  Next, we
 apply (\ref{37}) for $n_s^{(0)}$ and $n_0^{(1)}$ to obtain,
 with the help of  the identities (\ref{38}) - (\ref{100}),
\be
 \Omega(n_s^{(0)})=\lambda Q_0 + \mu Q_1
-an_s^{(0)}+bP, \quad \ein, \label{42} \ee

\nin  and \be
 \Omega(n_0^{(1)})=\lambda R_0 + \mu R_1
-bn_0^{(1)}+aP, \quad \ein. \label{43} \ee

\nin Now  we derive from (\ref{36}) the expression for $g_2$ that
we will be working with:

\be g_2(\eta)=-(\lambda+\mu)g_1(\eta)+\lambda aQ_0- \mu bR_1+ \mu
aQ_1-\lambda b R_0- a^2n_s^{(0)}+b^2n_0^{(1)},\quad \ein.
 \label{44}\ee

\nin (\ref{44}) and (\ref{36}) show that for the model (\ref{96})
the relationship  (\ref{90}) holds iff

\ber T(\eta):= \lambda aQ_0- \mu bR_1+ \mu aQ_1-\lambda b R_0-
a^2n_s^{(0)}+b^2n_0^{(1)}&-&\nonumber\\
(an_s^{(0)}-bn_0^{(1)})A_2- B_2\vert\eta\vert- C_2 =0,\quad
\ein,\label{45} \eer

 \nin where $A_2,B_2,C_2$ are  coefficients
that do not depend on $\ein.$ We put in (\ref{45}) first
$\eta=\emptyset$ and then $\eta=\bar{\emptyset}$ to obtain
$C_2=-\lambda bN$ and $B_2=\lambda b+\mu a.$

\nin We will treat separately the  case $ s\ge 3$ and the case
$s=2.$ Since $a=b=0$ leads to a particular  case of noisy voter
model, we suppose in the sequel that $a^2+b^2\neq 0.$ Let some
fixed $y,z\in V$ obey the condition (ii), and $u_1,u_2\in V$ are
defined as in (\ref{101}).

\nin The  case $s\ge 3$.
 Consider the following two configurations: $\eta_1,$ defined by
$\eta_1(y)=\eta_1(u_1)=0, \ \eta_1(v)=1,$  for all $v\neq u_1,y, $
and  $\eta_2=(\eta_1)_{u_2}.$
 It is easy to figure out the following relationships

\ber Q_0(\eta_i)=2,\quad i=1,2,\quad Q_1(\eta_1)=N-2s, \quad
Q_1(\eta_2)=N+1-3s,& & \nonumber\\ R_0(\eta_i)=R_1(\eta_i)=0,\quad
\quad n_s^{(0)}(\eta_i)=n_0^{(1)}(\eta_i)=0, \quad i=1,2. \la{48}
\eer

\nin Substituting (\ref{48}) in (\ref{45}) gives

\be
2\lambda a+\mu a(N-2s)-(N-2)(\lambda b+\mu a)+ \lambda bN=0 \ee

\nin and

\be
2\lambda a+\mu a(N+1-3s)-(N- 3)(\lambda b+\mu a)+ \lambda bN=0,
\ee

\nin which implies

\be
\mu a (s-1)=(a+b)\lambda, \quad a\mu(3s-4)=\lambda( 3b+ 2a).
\la{104A} \ee

\nin By the same argument, applied to the configurations
$\bar\eta_1, \bar\eta_2$ we also get

\be
\lambda b(s-1)=(a+b)\mu,\quad \lambda b(3s-4)=\mu (3a + 2b).
\la{105} \ee

\nin  We will find all  solutions of (\ref{104A}) and (\ref{105}).
First we see that $\lambda\mu \neq 0$ implies $a=b$ and
consequently, $\lambda=\mu>0,$  $s=3.$ This gives the threshold
voter model ($\mathbf{C_3}$). If $\lambda \mu=0,$ then we should
have $\lambda a=\lambda b=\mu a =\mu b=0.$ By (\ref{45}) this
implies $a=b$ or $ab=0.$  In the first case,  we have $\lambda
=\mu =0,$ which is again ($\mathbf{C_3}$), while in the second
case, $\lambda =\mu =ab=0,$ which is ($\mathbf{C_2}$).

\nin  The  case $s=2$. Taking $\eta_1$ as above, gives

\be Q_0(\eta_1)=2,\quad Q_1(\eta_1)=N-4, \quad R_1(\eta_1)=2,
\quad R_0(\eta_1)=n_2^{(0)}(\eta_1)=n_0^{(1)}(\eta_1)=0. \la{105A}
\ee

\nin Consequently, (\ref{45}) implies $(a+b)(\lambda-\mu)=0.$ If
$a+b=0$, then (\ref{45}) becomes

\ber \lambda a(Q_0+R_0)+ \mu
a(Q_1+R_1)-a^2(n_2^{(0)}-n_0^{(1)})&-&\nonumber\\
(n_2^{(0)}+n_0^{(1)})A_2a-a(\mu-\lambda)\vert\eta\vert-\lambda
aN=0,\quad \ein. \la{48a}\eer

\nin Since in the case $s=2$

\be Q_0+R_0=-2n_2^{(0)}-2n_0^{(1)}+N-\vert\eta\vert, \quad \ein
\ee

\nin and
\be
Q_1+R_1=-2n_2^{(0)}-2n_0^{(1)}+\vert\eta\vert, \quad \ein,
 \ee

\nin we see that (\ref{48a}) implies $a=0,$ and, consequently,
$b=0.$

\nin Let now  $\lambda=\mu.$ Then, we employ (\ref{R}), (\ref{39})
and (\ref{40}) to rewrite  (\ref{45})  as

\ber T(\eta)=\lambda\Big[aL(\eta)-bL(\bar\eta)\Big]
-a^2n_2^{(0)}+b^2n_0^{(1)}&-&\nonumber\\
(an_2^{(0)}-bn_0^{(1)})\tilde{A}_2=0, \quad \ein,\la{106}\eer

\nin where we denoted $L(\eta)=n_1^{(0)}+
n_2^{(1)}-\vert\eta\vert$ and $\tilde{A}_2=A_2+ 3\lambda.$ So,

\ber T(\eta)+T(\bar\eta)=(a-b)\Big[\lambda(L(\eta)+L(\bar\eta))
-(a+b)(n_2^{
(0)}+n_0^{(1)})&-&\nonumber\\
(n_2^{(0)}+n_0^{(1)})\tilde{A}_2\Big]=0, \quad \ein.\la{107} \eer

\nin First observe that  if $a=b$ and  $\lambda=\mu$ then we have
the model
 ($\mathbf{C_3}$). Next, substituting in (\ref{107})
$$L(\eta)+L(\bar\eta)=n_1^{(0)}+ n_2^{(1)}+ n_1^{(1)}+
n_0^{(0)}-N= -n_2^{(0)}- n_0^{(1)},$$
 \nin we have

\be
(n_2^{(0)}+ n_0^{(1)})(\tilde{A}_2+a+b+\lambda)=0, \quad \ein. \ee

\nin Consequently, $\tilde{A}_2=-a-b-\lambda,$ which in view of
(\ref{106}) yields

\be
\lambda a(n_1^{(0)}+n_2^{(1)}+n_2^{(0)} -\vert\eta\vert) -\lambda
b(n_1^{(1)}+n_0^{(0)} +n_0^{(1)}-N+\vert\eta\vert)
+ab(n_2^{(0)}-n_0^{(1)})=0, \quad \ein. \la{109}\ee

\nin A specific feature of the case $s=2$ is that the following
identity holds

\be
2(n_2^{(0)}-n_0^{(1)})= n_1^{(1)}-n_1^{(0)}, \quad \ein. \la{108}
\ee

\nin So, we obtain from (\ref{109})

\be
(n_2^{(0)}-n_0^{(1)})(-\lambda a -\lambda b +ab)=0, \quad \ein,
\ee

\nin which gives the model ($\mathbf{C_4}$).

 \nin This completes the proof of the necessity of the conditions
 $(\mathbf{C_1}) -(\mathbf{C_4}).$

\nin  The proof that each of the conditions $(\mathbf{C_1})
-(\mathbf{C_4})$ is sufficient for (\ref{cd})
    is now simple. In
the case ($\mathbf{C_1}$) it was shown in \cite{GR} that
\be
g_1(\eta)= \lambda_0 N -(\lambda_0+\mu_s)\vert\eta\vert, \quad
\ein, \ee

\nin which implies $g_2(\eta)= -(\lambda_0+\mu_s)g_1(\eta), \
\ein$ .

\nin In the case ($\mathbf{C_3}$) we have  $s=2,3,$
$\lambda_k=\mu_{k+1}=h, \quad k=0,\ldots,s-1, \quad
\lambda_s=\mu_0:=h+a,$ where $a\in R:h+a\ge 0.$ So, (\ref{44})
becomes

\be g_2(\eta)=-2hg_1(\eta)+ah( Q_0-  R_1+ Q_1- R_0)-
a^2(n_s^{(0)}-n_0^{(1)}),\quad \ein.
 \label{110}\ee

\nin  By the Lemma, (\ref{39}), (\ref{40}) and (\ref{R}) we
further obtain for $s=2,3$

\be
g_2(\eta)=-2h g_1(\eta)+ ha
\Big((F_1-s-1)(n_s^{(0)}-n_0^{(1)})+2\vert\eta\vert- N\Big)-
a^2(n_s^{(0)}-n_0^{(1)}),\quad \ein, \ee

\nin where $F_1$ is given by (\ref{89}). Hence, in view of
(\ref{36}), we have in  the both cases of $s$
\be
g_2(\eta)=-(a+8h)g_1-6h^2(2\vert\eta\vert-N),\quad \ein. \ee

\nin  Let now ($\mathbf{C_4}$) hold. With the help of (\ref{108})
it is easy to verify that (\ref{106}) indeed holds with
$\tilde{A}_2=-a-b-h.$ Finally, in the case of the model
($\mathbf{C_2}$) we have either $g_2(\eta)=-bg_1(\eta),\ \ein$ or
$g_2(\eta)=-ag_1(\eta),\ \ein.$ $\clubsuit$

\nin {\bf Corollary 1.} The mean coverage functions
$M^{(\eta)}(t),\quad \ein, \quad t\ge 0$ of the models
($\mathbf{C_i}$), $i=1,2,3, 4$ are given by the expressions
($\mathbf{D_i}$), $i=1,2,3,4$ correspondingly:

\nin ($\mathbf{D_1}$): \be
M^{(\eta)}(t)=(\vert\eta\vert-\frac{\lambda_0 N}{\lambda_0+\mu_s})
\exp{\Big(-(\lambda_0+\mu_s)t\Big)}+\frac{\lambda_0
N}{\lambda_0+\mu_s} , \quad t\ge 0, \quad \ein \ee

\nin ($\mathbf{D_2}$): \be
M^{(\eta)}(t)=-n_s^{(0)}(\eta)\exp{(-at)}+\vert\eta\vert +
n_s^{(0)}(\eta), \quad t\ge 0,\quad \ein\ee

\nin or \be
M^{(\eta)}(t)=n_0^{(1)}(\eta)\exp{(-bt)}+\vert\eta\vert -
n_0^{(1)}(\eta), \quad t\ge 0,\quad \ein.\ee

\nin ($\mathbf{D_3}$):
\be
M^{(\eta)}(t)= C_1 e^{-\alpha_1 t} +C_2 e^{-\alpha_2 t} +
\frac{N}{2}, \quad t\ge 0, \quad \ein, \ee

\nin where

\be
\alpha_i=\frac{8h +a -(-1)^i\sqrt{(8h+a)^2 - 48h^2}}{2}, \quad
i=1,2. \la{202}\ee

\nin and \be
C_1=\frac{\alpha_2(\vert\eta\vert-\frac{1}{2}N)+g_1(\eta)}
{\alpha_2-\alpha_1}, \quad
C_2=\vert\eta\vert-\frac{N}{2}-C_1,\quad \ein. \la{111}\ee

\nin ($\mathbf{D_4}$):

\be
M^{(\eta)}(t)= C_1 e^{-\alpha_1 t} +C_2 e^{-\alpha_2 t} +
 N\frac{4h+a}{8h+a+b},
\quad t\ge 0, \quad \ein, \ee

\nin where
\be
\alpha_i=\frac{6h+a+b -(-1)^i\sqrt{4h^2+ (a+b)^2 +8h(a+b)}}{2},
\quad i=1,2, \quad ab=h(a+b) \la{302}
 \ee

\nin and \be C_1=\frac{\alpha_2(\vert\eta\vert-
N\frac{4h+a}{8h+a+b})+g_1(\eta)}{\alpha_2-\alpha_1}, \quad
C_2=\vert\eta\vert- N\frac{4h+a}{8h+a+b}-C_1,\quad \ein.
\la{1111}\ee

\nin {\bf Proof.}

\nin The assertions  follow from the relationships
($\mathbf{E_i}$), $i=1.2,3,4$ below that hold for the models
 ($\mathbf{C_i}$), $i=1,2,3,4$ correspondingly.

\nin ($\mathbf{E_1}$): \be g_1(\eta)= \lambda_0N
-(\lambda_0+\mu_s)\vert\eta\vert, \quad \ein \ee

\nin ($\mathbf{E_2}$):
\be
g_2(\eta)=-ag_1(\eta), \quad a\ge 0, \quad \ein \ee

\nin or

\be g_2(\eta)=-bg_1(\eta), \quad b\ge 0, \quad\ein \ee

\nin ($\mathbf{E_3}$):
\be
g_2(\eta)=-(a+8h)g_1-6h^2(2\vert\eta\vert-N),\quad \ein \ee

\nin ($\mathbf{E_4}$):

\be
g_2(\eta)=-(6h+a+b)g_1(\eta)-h(8h+a+b)\vert\eta\vert+hN(4h+a),
\quad \ein. \ee

\nin Namely, the expressions ($\mathbf{D_1}$)-($\mathbf{D_4}$) are
obtained by solving the second order differential equations
corresponding to ($\mathbf{E_1}$)- ($\mathbf{E_4}$), under the
initial conditions
\be
M^{(\eta)}(0)=\vert\eta\vert, \ \
\frac{dM^{(\eta)}(0)}{dt}=g_1(\eta), \quad \ein. \ee
\nin \section{ The mean density function}

\nin Let $\nu$ be a probability measure on the state space
${\cal{X}}_N.$ Denote by $\varphi^{(\nu)}_t, \quad t\ge 0$  the
SNNSS starting from $\nu$ ( this means that  the distribution of
$\varphi^{(\nu)}_0$ is $\nu$), and denote by

\nin $M^{(\nu)}(t)=E_\nu M^{(\eta)}(t), \ t\ge 0$ the
corresponding mean coverage function. The function
$w_N^{(\nu)}(t)=N^{-1}M^{(\nu)}(t), \  t\ge 0$ is called the mean
density coverage function corresponding to the initial
distribution $\nu.$

\nin {\bf Historical remark.} The mean density function was
studied in a number of papers. In addition to the previously
mentioned literature that is immediately related to the context of
the present paper,  we outline now some  adjacent topics of
research. Special attention was devoted to the contact process.
 Gray \cite{G} investigated the behavior  of the population
profile function $p_t(x):=P(\varphi^{(\eta)}_t(x))=1, \ x\in V, \
t\ge 0,$ when $G=Z$ and $\eta$ is the empty configuration flipped
at the vertex $x=0.$ Belitsky \cite{B} treated a special case of
the previously mentioned adsorption-desorption process, when
$\mu_0>0, \ \mu_k=0, \ k=1,2$ and  $G=Z.$  It was proven in
\cite{B}, that the function $w_N^{(\nu_0)}(t), \ t\ge 0, $ where
$\nu_0$ is the measure concentrated on the empty configuration,
possesses a saddle point. This  extends the result of \cite{GW}.
Continuing the discussion in \cite{GW}, Belitsky \cite{B} relates
the  above phenomenon to the violation of the classical Langmuir
law, known in physical chemistry. Note, that from  \cite{B},  as
well as \cite{BG1}, one can see how complicated is  the structure
of the iterations $g_i, \ i\ge 1 $ of the generator of  the
process considered. This explains the difficulties in the study of
the transient behavior of functionals of contact process even in
the case $G=Z.$ The problem becomes much simpler in the framework
of the mean-field theory, that corresponds to the case when $G$ is
a complete graph. In this case, a SNNSS conforms to the
birth-death process (see Granovsky and Zeifman \cite{GZ}). The
limiting behavior, of the density process
$N^{-1}\vert\varphi_t\vert,$ as $N\to \infty,$ when $\varphi_t, \
t\ge 0$ is the basic contact process, was extensively studied in
the literature. For the most recent review of the topic see
Durrett \cite{Dur}.

\nin In conclusion, we mention two papers devoted to  voter
models. Cox \cite{C} derived the limit of the density process for
the basic voter model on the torus in $Z^d$, under an appropriate
time scaling.
 Mountford \cite{M} considered a class of one- dimensional
 multitype IPS,  that are
    featured by the following property of its generator
 $\Omega:$
\be
\sup_{n,\eta}\vert\Omega f_n(\varphi_t)\vert\le const, \la{co}\ee

\nin where $f_n(\eta)=\sum_{x\in Z:\vert x\vert\le n} \eta(x).$ He
proved that under condition (\ref{co}) the coverage process $
f_n(\eta_t), \ t\ge 0$  is a martingale plus a term that is
negligible as $n\to \infty.$ In this sense these models  can be
viewed as a  generalization of the basic voter model. It should be
noted that the condition (\ref{co}) fails for all SNNSS
($\mathbf{C_1}$)-($\mathbf{C_4}$), except only the case of the
basic voter model.

\nin An important particular case of $\nu$ is the product
Bernoulli measure $\nu_p,\ 0\le p\le 1, $ defined by

\be \nu_p(\eta)=\prod_{x\in V}p^{\eta(x)}(1-p)^{1-\eta(x)}=
p^{\vert\eta\vert}(1-p)^{N-\vert\eta\vert},\quad \ein.
 \ee

\nin For a given $0\le p\le 1,$ the measure $\nu_p$ corresponds to
the initial distribution on ${\cal{X}}_N,$ such that all spins are
i.i.d. Bernoulli random variables. In view of this,
$M^{(\nu_p)}(t)=NE\varphi^{(\nu_p)}_t(x), \ \forall x\in V,\ \ein
,$ and consequently, $w_N^{(\nu_p)}(t)=E\varphi^{(\nu_p)}_t(x), \
\forall x\in V, \ \ein.$ Hence, the mean density
$w_N^{(\nu_p)}(t), \ t\ge 0$ defines the marginal distribution of
the process $\varphi^{(\nu_p)}_t$ at any site $x\in V$ at time
$t\ge 0.$ It turns out that the densities $ w_N^{(\nu_p)}(t), \
t\ge 0$ corresponding to SNNSS's $(\mathbf{C_1}) -(\mathbf{C_4})$
, have the following remarkable property.

\nin {\bf Proposition 2.} The mean density functions
$w_N^{(\nu_p)}(t),\ t\ge 0, \  0\le p\le1$ of models
$(\mathbf{C_1}) -(\mathbf{C_4})$ do not depend on $N.$

\nin {\bf Proof.}

\nin It follows from (\ref{90}) that

\be
g_{i}(\eta)=A_{1}g_{i-1}(\eta) + A_{0}g_{i-2}(\eta), \quad
i=3,\ldots, \quad \ein, \la{120} \ee

\nin where, in view of ($\mathbf{E_1}$)- ($\mathbf{E_4}$), the
coefficients $A_{1},  \ A_{0}$ do not depend neither on $\ein$ nor
$N.$

\nin  We deduce from (\ref{05}) that

 \be
E^{(\nu_p)} g_1(\eta)=
NE^{(\nu_p)}\Big[c(x,\eta)(1-2\eta(x))\Big], \quad \forall x\in V.
\la{121}\ee

\nin It is clear that the expected value in the RHS of (\ref{121})
does not depend on $N,$ for any SNNSS. Next, in (\ref{90}) the
coefficient  $B=N B_0$, where, by ($\mathbf{E_1}$)-
($\mathbf{E_4}$), the factor $B_0$ does not depend on $N.$ So,
(\ref{90}) implies  $E^{(\nu_p)} g_2(\eta)= N q_2,$ where the
factor $q_2$ does not depend on $N.$ This together with
(\ref{120}) and (\ref{121}) gives $E^{(\nu_p)} g_i(\eta)= N q_i, \
i=3, \ldots,$ where again  the factors $q_i, \ i=3,\ldots,$ do not
depend on $N.$ Finally, to complete the proof, we use the
Hille-Yosida series expansion of the function $M^{(\nu_p)}(t), \
t\ge 0.$ $\clubsuit$

\nin {\bf Remark 1.} It can be shown that even linear SNNSS given
by $\lambda_k=\lambda_0+d_\lambda k, \ \ \mu_k=\mu_0 +d_\mu k, \
k=0,1,\ldots, s$ with $d_\lambda\neq d_\mu $  have not the
property stated in Proposition 2.

\vskip .4cm
 \nin In view of  Proposition 2,  we write
$w^{(\nu_p)}(t):=w_N^{(\nu_p)}(t), \ t\ge 0, \ 0\le p\le 1, \
N=s+1,\ldots,  $  for the models $(\mathbf{C_1}) -(\mathbf{C_4})$.
The explicit expressions for the functions $w^{(\nu_p)}(t), \ t\ge
0$ is easy to obtain from ($\mathbf{E_1}$)- ($\mathbf{E_4}$). Now
the Trotter - Kurtz approximation theorem (\cite{Lig1}) and
($\mathbf{E_1}$)- ($\mathbf{E_4}$) give immediately the following

\nin {\bf Corollary 2.} For models $(\mathbf{C_1})
-(\mathbf{C_4})$ on a finite or infinite $s$-regular graph,
 \be
Pr(\varphi_t^{(\nu_p)}(x)=1)=w^{(\nu_p)}(t),\quad t\ge 0, \quad
x\in V, \quad 0\le p\le 1. \ee

\nin In particular, in the presence of noise, the processes
($\mathbf{C_1}$),($\mathbf{C_3}$), ($\mathbf{C_4}$) have  ergodic
marginals, in the sense that for each of these processes

\be
\lim_{t\to \infty} w^{(\nu_p)}(t) \la{li} \ee

\nin exists and does not depend on  $0\le p\le 1.$

\nin  The expressions ($\mathbf{D_1}$),($\mathbf{D_3}$) and
($\mathbf{D_4}$) give  correspondingly the following values for
the limit in (\ref{li}):

\be \frac{\lambda_0}{\lambda_0+\mu_s}, \quad \frac{1}{2} \quad
and\quad \frac{4h+a}{8h+a+b}.\ee

\nin {\bf Remark 2.} Formulae ($\mathbf{D_1}$)-($\mathbf{D_4}$)
show the complicated influence of a constant additive  noise  on a
transient behaviour of the process. In particular, note that under
the absence of noise

\nin ( $h=0$) we have in (\ref{202}), (\ref{302}) $\alpha_2=0.$ It
is also appropriate to mention that the processes
 ($\mathbf{C_1}$)- ($\mathbf{C_4}$) are either attractive
or anti-attractive. The latter means that in the definition of
attractiveness ( \cite{Lig1}, p.132) the direction of inequalities
for the
 flip rates is  reversed.
\nin \section{Ergodicity and Spectral gap}

\nin {\bf Ergodicity.} We will address the question of ergodicity
of the processes ($\mathbf{C_1}$)-($\mathbf{C_4}$) on a finite or
infinite $s$-regular graph. The process ($\mathbf{C_2}$) is,
obviously, not ergodic. The processes ($\mathbf{C_1}$),
($\mathbf{C_3}$) with $a\ge 0$ and ($\mathbf{C_4}$) with $a,b\ge
0$  are  attractive.  The key property of such processes is
 that ergodicity of their  marginals implies the ergodicity of the
 process. A beautiful argument leading to this assertion is
 explained
 in \cite{Lig1}(see Corollary 2.8, p. 75, and Corollary 2.4, p.136).
\nin So, by Corollary 2, we get

\nin {\bf Corollary 3.} The following three processes are ergodic:
($\mathbf{C_1}$) with $h_1+h_2=\lambda_0+\mu_s>0,$
($\mathbf{C_3}$) with $a\ge 0, h>0$ and ($\mathbf{C_4}$) with
$a,b\ge 0,\ h>0.$

\nin {\bf Remark 3} i.The ergodicity of the first among the three
models in Corollary 3 was proven in \cite{GR}. The ergodicity of
the third one as well as  the second one in the case  $s=2,$
follows from the fact that these are attractive  spin systems in
one dimension with translation invariant and positive flip rates
(see \cite{Lig1}, Theorem 3.14, p.152). To the best of our
knowledge, the established ergodicity of the  second model in the
case $s=3$ answers an open question. We will explain below that
the $\epsilon -M>0$ condition in the case considered gives
ergodicity for $a<\frac{2}{3}h$ only.

\nin ii. It is interesting to observe that, by Corollary 2, the
models considered have ergodic marginals also in the case when
they are not attractive and even not ergodic. For example, this is
true for ($\mathbf{C_3}$), when  $h+a=0,\ s=2,$ in which case the
process on $Z_1$ is not ergodic, having two different absorbing
states $\eta_i: \eta_1(x)=0.5(1+(-1)^{i+\vert x\vert}), \ i=1,2, \
x\in Z_1. $

\nin {\bf Spectral Gap.} Recall (see for references  \cite{CH1},
\cite{CH2}) that the spectral gap $\alpha>0 $ of an exponentially
ergodic Feller-Markov process $\varphi_t\quad t\ge 0$ with an
invariant measure $\nu$ on  state space $\cal X=\{\eta\}$ is
defined by
\be
\alpha=\sup\{\beta>0:\sup_{\eta\in {\cal X}}\vert
Ef(\varphi_t^{(\eta)})- \int f d\nu\vert\le A_f\exp(-\beta t),
\quad f\in {C(\cal X)} , \quad t\ge 0 \}, \la{200} \ee

\nin where $A_f\ge 0 $  does not depend on $t\ge 0.$

\nin It is plain that formulae  ($\mathbf{D_1}$), ($\mathbf{D_3}$)
and ($\mathbf{D_4}$) for the mean coverage function provide an
estimate from above   for  spectral gaps of the corresponding
processes. Namely, if $\tilde{\alpha}$ is the rate of exponential
convergence to the equilibrium of the function $M^{(\eta)}(t), \
t\ge 0$ (=$w^{(\nu_p)}(t), \ t\ge0   )$, then
$\alpha\le\tilde{\alpha}.$ From the other hand, the celebrated
$\epsilon-M>0$ condition of ergodicity (\cite{Lig1}, p.31)
provides the lower bound of the spectral gap of any ergodic IPS on
a finite or infinite graph: $\alpha\ge\epsilon-M>0.$ In the case
of a SNNSS on a $s$-regular graph the quantities $\epsilon$ and
$M$ are given by
\be
\epsilon= \min_{0\le k\le s}\left\{\lambda_k+\mu_k\right\} \quad
M=s
\max_{0\le k\le s}
\left\{\vert\lambda_k-\lambda_{k-1}\vert,\vert\mu_k-
\mu_{k-1}\vert\right\}.\ee

 \nin Below are the values of
$\epsilon -M$ and  $\tilde{ \alpha}$ for the three ergodic models
considered.

\nin ($\mathbf{C_1}$):
\be
\epsilon -M=\lambda_0+\mu_s=\tilde{\alpha} \ee

\nin So,   \be\alpha=\lambda_0+\mu_s, \la{201}\ee

\nin if $\lambda_0+\mu_s>0.$ This fact was observed in \cite{GR}.
Moreover, it was proven in \cite{GZ} that if $G$ is a complete
graph, then the noisy voter model is the only one SNNSS with the
property (\ref{201}).

\nin ($\mathbf{C_3}$): \be\tilde{\alpha}=\alpha_1,\ee

\nin where $\alpha_1$ is given by (\ref{202}), while

\be \epsilon -M=\left\{\begin{array}{ll} 2h-sa, & {\rm if~} a>0
\cr 2h+(s+1)a, & {\rm if~} a\le 0,
\end{array}
\right.
\ee

\nin where $s=2,3$. This says that in the case $a>0$ the $\epsilon
-M>0$ condition is  applicable  for $s=2,$  if $h>a,$ and for
$s=3,$ if $h>1.5 a.$  It is easy to verify that $\alpha_1=\epsilon
-M>0$ iff $a=0$ which is the trivial case of the noisy voter
model.

 \nin ($\mathbf{C_4}$):
\be\tilde{\alpha}=\alpha_1,\ee

\nin where $\alpha_1$ is given by (\ref{302}), while

\be \epsilon -M= 2h +\min\{0,a,b\}-2\{\max\vert a\vert,\vert
b\vert\}, \la{301}\ee

\nin where $h=\frac{ab}{a+b}>0.$

\nin Here again one can see that $\alpha_1=\epsilon -M>0$ is
impossible unless $a=b=0.$

\vskip .5cm

 \nin The preceding discussion leads us to the
following

\nin {\bf Conjecture.} Noisy Voter model is the only one NNSS for
 which $\alpha=\epsilon -M>0.$
\vskip .5cm

\nin {\bf Concluding remark.
} It is demonstrated by our Theorem,
that passing from the first order differential equation to the one
of the second order does not enrich much the class of

\nin solvable ( in the sense of the mean coverage function )
SNNSS. A natural question arising in this connection is the
characterization of SNNSS in the case of higher order differential
equations. One can expect that the progress in this direction will
lead to the discovery of wider classes of solvable SNNSS. The
solution to this problem requires the analysis of the structure of
generators $g_i,$ as defined in (\ref{ab}), of higher orders.
 At the moment the problem looks intractable.

\vskip .4cm

\nin {\bf Acknowledgement} The research was supported by the fund
of the Promotion of Research at Technion and by the V.P. R. fund
at the Technion.

The author is thankful to Prof. Tom Liggett for his precious
comments on the first draft of the paper and to the two referees
for their constructive criticism that helped to improve the
exposition.

\end{document}